\documentclass[12pt,reqno]{amsart}
\usepackage{amsmath,amssymb}
\usepackage[dvips]{graphicx}
\usepackage[colorlinks,linkcolor=blue,citecolor=blue]{hyperref}
\newtheorem{thm}{Theorem}
\newtheorem{lem}[thm]{Lemma}
\newtheorem{exmp}[thm]{Example}

\newcommand{\mn}{\mbox{Min}}
\newcommand{\irr}{\mbox{Irr}}
\newcommand{\lcm}{\mbox{Lcm}}

\newcommand{\B}{\mbox{B}}

\newcommand{\coeff}{\mbox{Coeff}}

\newcommand{\merge}{\mbox{Merge}}
\newcommand{\minmerge}{\mbox{MinMerge}}
\newcommand{\maxmerge}{\mbox{MaxMerge}}
\newcommand{\rmv}[1]{}

\newcommand{\N}{\mathbb{N}}
\newcommand{\K}{\mathbb{K}}
\newcommand{\bN}{\overline{\mathbb{N}}}
\newcommand{\cT}{\mathcal{T}}

\begin{document}

\title{Computing Irreducible Decomposition of Monomial Ideals}

\thanks{The authors were partially
supported by the National Science Foundation under grant DMS-0302549
and National Security Agency under grant H98230-08-1-0030.}

\author {Shuhong Gao and Mingfu Zhu}
\address{Department of Mathematical Sciences, Clemson University, Clemson, SC 29634-0975 USA \\\{sgao,
mzhu\}@clemson.edu}

\begin{abstract}
The paper presents two algorithms for finding irreducible
decomposition of monomial ideals. The first one is recursive,
derived from staircase structures of monomial ideals. This algorithm
has a good performance for highly non-generic monomial ideals. The
second one is an incremental algorithm, which computes
decompositions of ideals by adding one generator at a time. Our
analysis shows that the second algorithm is more efficient than the
first one for generic monomial ideals. Furthermore, the time
complexity of the second algorithm is at most $O(n^2p\ell )$ where
$n$ is the number of variables, $p$ is the number of minimal
generators and $\ell$ is the number of irreducible components.
Another novelty of the second algorithm is that, for generic
monomial ideals, the intermediate storage is always bounded by the
final output size which may be exponential in the input size.
\end{abstract}


\maketitle
% Introduction ==============================================================
\section{Introduction}

Monomial ideals provide ubiquitous links between combinatorics and
commutative algebra \cite{Vil01,MS05}. Though simple they carry
plentiful algebraic and geometric information of general ideals. Our
interest in monomial ideals is motivated by a paper of \cite{GRS03},
where they studied the connection between the structure of monomial
basis and the geometric structure of the solution sets of
zero-dimensional polynomial ideals. Irreducible decomposition of
monomial ideals is a basic computational problem and it finds
applications in several areas, ranging from pure mathematics to
computational biology, see for example \cite{HS04} for computing
integer programming gaps, \cite{BY06} for computing tropical convex
hulls, \cite{SS06} for finding the joins and secant varieties of
monomial ideals, \cite{Anw07} for partition of a simplicial complex,
\cite{Rou08a} for solving the Frobenius problem, and \cite{JLS06}
for modeling gene networks.

We are interested in efficient algorithms for computing irreducible
decomposition of monomial ideals. There are a variety of algorithms
available in the literature. The so-called splitting algorithm:
Algorithm 3.1.2 in \cite{Vas98} is not efficient on large scale
monomial ideals. \cite{Mil04} gives two algorithms: one is based on
Alexander duality \cite{Mil00}, and the other is based on Scarf
complex \cite{BPS98}. \cite{Rou07} improves the Scarf complex method
by a factor of up to more than 1000. Recently, \cite{Rou08b}
proposed several slicing algorithms based on various strategies.

Our goals in this paper are to study the structure of monomial
ideals and present two new algorithms for irreducible decomposition.
We first observe some staircase structural properties of monomial
bases in Section \ref{structure}. The recursive algorithm presented
in Section \ref{recursive} is based on these properties, which allow
decomposition of monomial ideals recursively from lower to higher
dimensions. This algorithm was presented as posters in ISSAC 2005
and in the workshop on Algorithms in Algebraic Geometry at IMA in
2006. Our algorithm was recently generalized by \cite{Rou08b} where
several cutting strategies were developed and our algorithm
corresponds to the minimum strategy there. Also, the computational
experiments there shows that our algorithm has good performance for
most cases, especially for highly non-generic monomial ideals.

Our second algorithm is presented in Section \ref{distrbn}.
 It can be viewed as an improved Alexander dual method (\cite{Mil00,Mil04}).
It is incremental based on some distribution rules for
\lq\lq$+$\rq\rq\ and \lq\lq$\cap$\rq\rq\ operations of monomial
ideals. We maintain an output list of irreducible components, and at
each step we add one generator and update the output list. In
\cite{Mil04}, there is no specific criterion for selecting
candidates that need to be updated, and the updating process is
inefficient too. Our algorithm avoids these two deficiencies. Our
analysis in Section \ref{Conclusion} shows that the second algorithm
works more efficiently than the first algorithm for generic monomial
ideals. We prove that, for generic monomial ideals, the intermediate
storage size (ie.\ number of irreducible components at each stage)
is always bounded by the final output size, provided that the
generators are added in lex order. This enables us to show that the
time complexity of the second algorithm is at most $O(n^2p\ell)$
where $n$ is the number of variables, $p$ is the number of minimal
generators and $\ell$ is the number of irreducible components.

In Section \ref{prelim}, we present some notations and introductory
materials on monomial ideals. In Section \ref{treer} we discuss tree
representations and operations of monomial ideals.\\

% Prelim ================================================================
\section{Monomial Ideals}\label{prelim}

We refer the reader to the books of \cite{CLO97} for background in
algebraic geometry and commutative algebra, and to the monograph
\cite{MS05} for monomial ideals and their combinatorial properties.

Let $\K$ be a field and $\K[X]$, the polynomial ring over $\K$ in
$n$ indeterminates $X=x_1,\ldots, x_n$. For a vector
$\alpha=(a_1,\ldots,a_n) \in \N^n$, where $\N=\{0,1,2,\ldots \}$
denotes the set of nonnegative integers, we set
$$
X^\alpha=x_1^{a_1}\ldots x_n^{a_n},
$$
which is called a {\bf monomial}. Thus monomials in $n$ variables are in $1-1$ correspondence with vectors in
$\N^n$. Suppose $\alpha=(a_1,\ldots,a_n)$ and $\beta=(b_1,\ldots,b_n)$ are two vectors in $\N^n$,
we say
$$\alpha \le \beta \textrm{ if } a_j \le b_j \textrm{ for all } 1\le j \le n.$$ This defines a partial order on
$\N^n$, which corresponds to division order for monomials since $x^{\alpha}|x^{\beta}$ if and only if
$\alpha \le \beta$. We say
$$\alpha < \beta \textrm{ if } \alpha \le \beta \textrm{ but } \alpha \neq \beta.$$
Also we define
$$\alpha \prec \beta \textrm{ if } a_j <b_j \textrm{ for all } 1 \le j \le n.$$
Then $\alpha \nprec \beta$ means that $a_j \ge b_j$ for at least one
$j$.

An ideal $I \subset \K[X]$ is called a {\bf monomial ideal} if it is
generated by monomials. Dickson's Lemma states that every monomial
ideal in $\K[X]$ has a unique minimal set of monomial generators,
and this set is finite. Denote this set to be $\mn(I)$, that is,
$$\mn(I)=\{X^\alpha \in I: \textrm{ there is no }
X^\beta \in I\ \textrm{ such that } \beta < \alpha \}.$$

A monomial ideal $I$ is called {\bf Artinian} if $I$ contains a
power of each variable, or equivalently, if the quotient ring
$\K[X]/I$ has finite dimension as vector space over $\K$. For
convenience of notations, we define
\[ x_i^{\infty}=0, \quad 1 \le i \le n.\]
By adding infinity power of variables if necessary, a non-Artinian
monomial ideal can be treated like an Artinian monomial ideal. For
example, $I=\langle x^2y^3 \rangle=\langle x^{\infty}, x^2y^3,
y^{\infty} \rangle$. Instead of adding infinity powers, we can also
add powers $x_i^{c_i}$ where $c_i$ is a sufficiently large integer,
say larger than the largest degree of $x_i$ in all the monomials in
$\mn(I)$. Then the irreducible components of the original ideal are
in 1-1 correspondence to those of the modified Artinian ideal; See
Exercise 5.8 in \cite{MS05} or Proposition 3 in \cite{Rou08b}. In
our algorithms belows, we will use infinity powers, but in the
proofs of all the results, we will use powers $x_i^{c_i}$.

An ideal $J \subset \K[X]$ is called {\bf irreducible} if it can not
be expressed as the intersection of two strictly larger ideals in
$\K[X]$. That is, $J=J_1 \cap J_2$ implies that $J=J_1$ or $J=J_2$.
A monomial ideal $I$ is irreducible if and only if $I$ is of the
form
$$
m^{\beta}=\langle x_1^{b_{1}},\ldots,x_n^{b_{n}}\rangle
$$
for some vector $\beta=(b_{1},\ldots,b_{n}) \in \bN^n$ where $\bN=\N
\cup \{ \infty \} \setminus \{ 0 \}$. Thus irreducible monomial
ideals are in 1-1 correspondence with $\beta \in \bN^n$.

An {\bf irreducible decomposition} of a monomial ideal $I$ is an expression of the form
\begin{equation} \label{eqn1.1}
I=m^{\beta_1} \cap \cdots \cap m^{\beta_r}
\end{equation} where $\beta_1,\ldots,\beta_r \in \bN^n$.
Since the polynomial ring $\K[X]$ is Noetherian, every ideal can be
written as irredundant intersection of irreducible ideals. Such an
intersection is not unique for a general ideal, but unique for a
monomial ideal. We say that the irreducible decomposition
(\ref{eqn1.1}) is \textbf{irredundant} if none of the components can
be dropped from the right hand side. If (\ref{eqn1.1}) is
irredundant, then the ideals $m^{\beta_1},\ldots,m^{\beta_r}$ are
called {\bf irreducible components} of $I$. We denote by ${\irr(I)}$
the set of exponents of irreducible components of $I$, that is,
$$\irr(I)=\{\beta_1,\ldots,\beta_r\}.$$
By this notation, we have
\[ I = \bigcap_{\beta \in \irr(I)} m^{\beta}.\]

Note that, for two vectors $\alpha$ and $\beta$, $$X^\alpha \in
m^{\beta} \textrm{ if and only if } \alpha \nprec \beta,$$ and
$$m^\alpha \subset m^{\beta} \textrm{ if and only if } \beta \le
\alpha.$$

A monomial ideal $I$ is called \textbf{generic} if no variable $x_i$
appears with the same non-zero exponent in two distinct minimal
generators of $I$. This definition comes from \cite{BPS98}. For
example,
$$I_1=\langle x^4,y^4,x^3y^2z,xy^3z^2,x^2yz^3 \rangle$$
is generic, but $$I_2=\langle x^4,y^4,x^3y^2z^2,xy^3z^2,x^2yz^3
\rangle$$ is non-generic, as $z^{2}$ appears in two generators.
Loosely speaking, we can say $I_2$ is nearly generic, but
$$I_3=\langle xy,yz,xz, z^2\rangle$$ is highly non-generic. Previous
algorithms \cite{Mil04,Rou07} behave very different for generic
monomial ideals and highly non-generic monomial ideals. For example,
the Scarf complex method works more efficient when dealing with
generic monomial ideals \cite{Mil04}.

In the following sections, we always assume that we are given the
minimal generating set of a monomial ideal. Though our algorithms
work for monomial ideals given by an arbitrary set of generators, it
will be more efficient if the generators are made minimal first.\\

% tree ================================================================
\section{Tree Representation and Operations}\label{treer}
Note that monomials are represented by vectors in $\N^n$ and
irreducible components are represented by vectors in $\bN^n$. To
efficiently represent a collect of vectors, we use a tree structure.
This is used in \cite{GRS03,Mil04}. This
data structure is also widely used in computer science, where it is called a trie.\\

\noindent \textbf{Tree representation.} First we want to define the
orderings on $\N^n$ or $\bN^n$. Suppose $\alpha=(a_1,\ldots,a_n)$
and $\beta=(b_1,\ldots,b_n)$ are two vectors in $\N^n$ or $\bN^n$,
and the variable ordering is $x_1< \cdots <x_n$ in $\K[X]$. We say
$\alpha <_{lex} \beta$ if $a_j=b_j$ for $k+1 \le j \le n$, but
$a_{k}<b_{k}$ for some $1 \le k \le n$.

\begin{figure}[h]
$$
\includegraphics*[222,697][390,813]{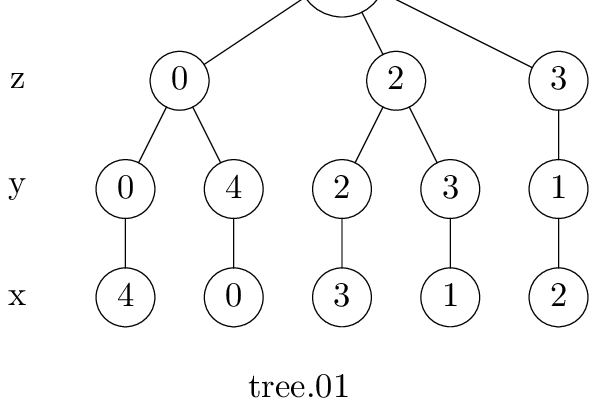}
$$
\caption{An example of tree representation.} \label{tree}
\end{figure}

Next, suppose $S \subset \N^n$ is a set of vectors corresponding to
the generators of a monomial ideal $I \subset \K[X]$. We represent
$S$ as a rooted tree $\cT$ of height $n$ in a natural way. The tree
should have $|S|$ leaves and the unique path of the tree from the
root to a leaf represents a vector in $S$. Precisely, to represent a
vector $\alpha=(a_1, \ldots, a_n)$, we label all the nodes except
the root of the path simply by $a_n, \ldots, a_1$ in the order from
the root to the leaf. We regard the root as being at height $0$. For
two vectors $\alpha=(a_1, \ldots, a_n)$ and $\beta=(b_1, \ldots,
b_n)$, if $a_j=b_j$ for $k+1 \le j \le n$ but $a_k \ne b_k$, then
$\alpha$ and $\beta$ share their corresponding path until height
$n-k$. After that their children are listed in increasing order with
respect to their coordinates. Figure \ref{tree} is the tree
representation for $I=\langle x^4,y^4,x^3y^2z^2,xy^3z^2,x^2yz^3
\rangle$ with variable order $x<y<z$.

The tree representation for a set of irreducible components could be constructed in a similar manner.

To perform the operations on sets of vectors, we need only perform on trees. We need three basic tree operations:
$\merge, \minmerge$
and $\maxmerge$.\\

\noindent \textbf{Merge.} Given $q$ rooted trees $\cT_1, \ldots, \cT_q$ with the same height,
merge them to form one rooted tree with the same height. Here we simply put the paths from all the trees together
with repetition ignored (actually no repeated paths occur in our algorithms). We
stress that no reduction work is performed under this operation.\\

\noindent \textbf{MinMerge.} We use $\minmerge(\cT_1, \ldots,
\cT_q)$ to represent the set of minimal elements in $\merge(\cT_1,
\ldots, \cT_q)$. For two vectors $\alpha, \beta$ in $\merge(\cT_1,
\ldots, \cT_q)$, if $\alpha \le \beta$, ie. $x^{\alpha}|x^{\beta}$,
then the path for $\beta$ should be removed in this operation. The
purpose is to find the minimal generating set for the ideal $I_{1} +
\cdots + I_{q}$ where
$\cT_{i}$ is the tree representation for $I_{i}$.\\

\noindent \textbf{MaxMerge.} Similarly, the set of maximal elements
in $\merge(\cT_1, \ldots, \cT_q)$ is represented by
$\maxmerge(\cT_1, \ldots, \cT_q)$. If $\alpha \le \beta$, ie.
$m^{\beta} \subset m^{\alpha}$, then the path for $\alpha$ should be
removed in this operation. Hence, if $\cT_{i}$ represents the set of
irreducible components of $I_{i}$, $1\leq i \leq q$, then
$\maxmerge(\cT_1, \ldots, \cT_q)$ represents the
 the set of irreducible components of the ideal $ I_{1} \cap \cdots \cap I_{q}$.\\

\section{Structure Properties of Monomial Bases} \label{structure}
In the results and their proofs below, we explicitly assume that all
the ideals are Artinian, adding large powers $x_i^N$ if necessary
where $N$ is an integer, though infinity powers will be used in the
Algorithms and Examples.

The monomial basis $\B(I)$ for a monomial ideal $I$ is defined as
$$\B(I)=\{\gamma\in \N^n: X^\gamma \notin I \},$$
which form a linear basis for the quotient ring $\K[X]/I$ over $\K$.
Thus, for $\gamma \in \N^{n}$, $\gamma \in B(I)$ if and only if
$\alpha \nleq \gamma$ for every $\alpha \in \mn(I)$. Note that
$\B(I)$ is a $\delta$-set, that is, if $\gamma \in \B(I)$ and $\mu
\le \gamma$, then $\mu \in \B(I)$. The next lemma characterizes
$B(I)$ in terms of $\irr(I)$.
\begin{lem} \label{bi}
For $\gamma\in \N^n$, $\gamma \in B(I)$ if and only if
$\gamma \prec \beta$ for some $\beta \in \irr(I)$.
\end{lem}
\begin{proof} Since $I = \cap_{\beta \in \irr(I)} m^{\beta}$, we have
 $X^{\gamma} \in I$ if and only if $X^{\gamma} \in m^{\beta}$, ie., $\gamma \nprec \beta$,
 for each $\beta \in \irr(I)$. Hence $X^{\gamma} \notin I$ if and only if
 $\gamma \prec \beta$ for some $\beta \in \irr(I)$, as desired.
\end{proof}

We now want to express $\irr(I)$ in terms of $\B(I)$. Since $I$ is
Artinian, for $\beta =(b_1,\ldots,b_n)\in \irr(I)$, we have $b_i>0$
for $1 \leq i \leq n$. Define
$$\beta \ominus 1=(b_1-1,b_2-1,\ldots,b_n-1).$$
Lemma \ref{bi} implies that, for each $\beta\in \irr(I)$, we have
$\beta\ominus 1 \in \B(I)$.

A vector $\gamma \in \N^{n}$ is called {\it maximal} in $\B(I)$ if
$$\gamma \in \B(I) \textrm{ and there is no } \mu \in \B(I) \textrm{ such that } \mu>\gamma.$$

\begin{lem} \label{irri}
For any vector $\beta \in \N^{n}$, $\beta\in \irr(I)$ if and only if
$\beta \ominus 1$ is maximal in $\B(I)$.
\end{lem}
\begin{proof} By Lemma \ref{bi}, $\beta \ominus 1 \in \B(I)$ if and only if there is
$\alpha \in \irr(I)$ such that $\beta \ominus 1 \prec \alpha$.
Notice that $ \alpha\ominus 1 \in \B(I)$ and $\beta \ominus 1 \prec
\alpha$ is equivalent to say $\beta \ominus 1 \le \alpha\ominus 1$.
Hence $\beta \ominus 1$ is maximal in $B(I)$ if and only if $\beta
\ominus 1=\alpha \ominus 1$, that is,
 $\beta=\alpha\in \irr(I)$.
\end{proof}

\begin{figure}[h]
$$
\scalebox{0.9}[0.9]{\includegraphics*[210,633][400,813]{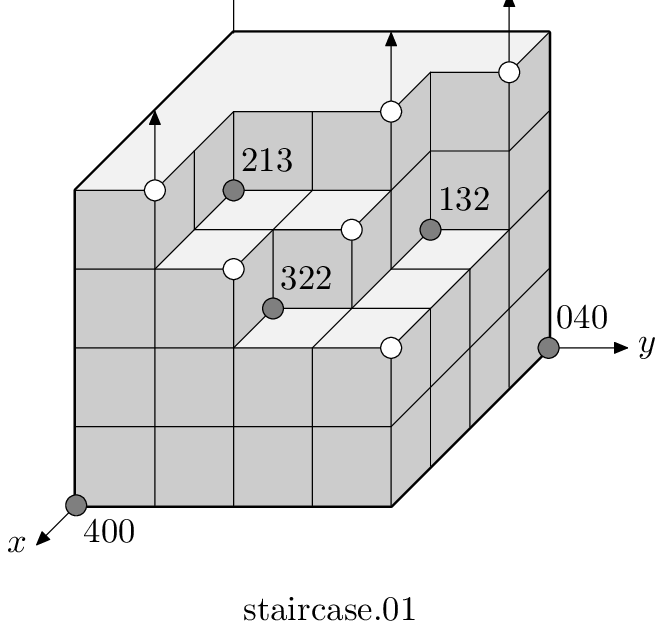}}
$$
\caption{An example of staircase diagram.} \label{irr}
\end{figure}

The staircase diagram will help us visualize the structural
properties of monomial ideals. For example, Figure \ref{irr} is the
staircase diagram for the monomial ideal $I=\langle
x^4,y^4,x^3y^2z^2,xy^3z^2,x^2yz^3 \rangle$. In this figure the gray
points are in 1-1 correspondence with the minimal generators, while
the white points are in 1-1 correspondence with the irreducible
components of $I$. Geometrically, $\B(I)$ is exactly the set of
interior integral points of the solid.

\section{Recursive Algorithm} \label{recursive}
For bivariate monomial ideals, irreducible decomposition is simple
\cite{MS99}. Suppose
$$\mn(I) =\left\{ x^{a_1}, x^{a_{2}}y^{b_2},\ldots,x^{a_{p-1}}y^{b_{p-1}}, y^{b_p}\right\}$$
where $a_1 > \cdots > a_{p-1} >0$, $0<b_2 < \cdots < b_p$, and $a_1$ or $b_p$ can be infinity. Then the irreducible
decomposition of $I$ is
$$I=\langle x^{a_1},
y^{b_2}\rangle \cap\langle x^{a_{2}}, y^{b_3}\rangle\cap \cdots \cap \langle x^{a_{p-2}}, y^{b_{p-1}}\rangle \cap
\langle x^{a_{p-1}}, y^{b_p}\rangle.$$

Our recursive algorithm is a generalization of the above observation
to higher dimensions. Let $I \subset \K[x_1,\ldots, x_{n}]$ be a
monomial ideal. Suppose all the distinct degrees of $x_n$ in
$\mn(I)$ are
$$0=d_0 <d_1< \ldots <d_s.$$ For example, in $I=\langle x^2y^3 \rangle=\langle x^{\infty}, x^2y^3, y^{\infty}
\rangle$, the distinct degrees in $y$ are $d_0=0,d_1=3$ and $d_3=\infty$. We collect the coefficients of $m \in
\mn(I)$ as polynomials in $x_{n}$. Precisely, for $0 \le k \le s$, let
$$I_k=\langle \coeff_{x_n}(m):\ m \in \mn(I) \textrm{ and } \deg_{x_n}
m \le d_k\rangle \subseteq \K[x_{1}, \ldots, x_{n-1}].$$ Then
\begin{equation} \label{ichn} I_0 \subsetneq I_1 \subsetneq \cdots \subsetneq I_{s}.
\end{equation}
By (\ref{ichn}), it follows that $$B(I_0) \supsetneq B(I_1)
\supsetneq \cdots \supsetneq B(I_s).$$ For the example with
$I=\langle x^{\infty}, x^2y^3, y^{\infty} \rangle$, $I_0=\langle
x^{\infty} \rangle=\{0\}$, $I_1=\langle x^{\infty}, x^2
\rangle=\langle x^2 \rangle$, and $I_2=\langle x^{\infty}, x^2, 1
\rangle=\langle 1 \rangle=\K[x]$.

We show how to read off the irreducible components of $I$ from those
of $I_{k}$'s, which have one less variables. For any vector
$\mu=(u_1,\ldots,u_{n-1})\in \N^{n-1}$ and $d\in \N$, define
$$(\mu,d)=(u_1,\ldots,u_{n-1}, d)\in \N^{n}.$$
\begin{lem}\label{d}
For any $\mu \in \N^{n-1}$ and $d \in \N$,
 $(\mu,d) \in B(I)$ if and only if there exists $k$, where $1 \le k \le s$, such that $d_{k-1} \le d< d_k$
 and $\mu \in B(I_{k-1})$.
\end{lem}
\begin{proof} $(\mu,d) \in B(I)$ if and only if there is no $m\in \mn(I)$ such that $m|X^{(\mu,d)}$. As $d_{k-1} \le d<
d_k$, we only need to see that there is no $m\in \mn(I)$ with
$\deg_{x_n} m \le d_{k-1}$. But this is equivalent to requiring that
$\mu \in B(I_{k-1})$.
\end{proof}

For a set of vectors $U$ and an integer $d$, define
$$U \otimes d=\{(u,d): u\in U\}.$$
\begin{thm} \label{rec} $\irr(I)
=\bigcup_{k=1}^{s}\big(\irr(I_{k-1})\setminus \irr(I_k)\big)\otimes d_k$, which is a disjoint union.
\end{thm}
\begin{proof} Assume $\mu \in \irr(I_{k-1})\setminus \irr(I_k)$. We first show that
 $(\mu,d_k) \ominus 1 \in B(I)$ and $\mu \ominus 1 \in B(I_{k-1})\setminus B(I_k)$.
Since $\mu \in \irr(I_{k-1})$, we have $\mu \ominus 1 \in
B(I_{k-1})$, so $(\mu,d_k) \ominus 1 = (\mu \ominus 1, d_k-1) \in
\B(I)$ by Lemma \ref{d}. Also, by Lemma \ref{irri}, there is no
$\gamma \in B(I_{k-1})$ such that $\gamma > \mu \ominus 1$, in
particular no $\gamma \in B(I_k)$ such that $\gamma > \mu \ominus
1$, as $B(I_k) \subset B(I_{k-1})$. Thus $\mu \ominus 1 \notin
B(I_k)$, otherwise we would have $\mu \in \irr(I_k)$ which
contradicts the assumption on $\mu$.

For $(\mu,d_k) \in \irr(I)$, we need to prove that $(\mu,
d_k)\ominus 1$ is maximal in $B(I)$. Assume otherwise, say $(\gamma,
d)\in B(I)$ and $(\gamma, d)>(\mu, d_k)\ominus 1$. Then $d \ge d_k$
or $d=d_k-1$. If $d \ge d_k$, then $\gamma \in B(I_j)$ where $k \le
j \le s$ by Lemma \ref{d}. Since $\gamma \ge \mu \ominus 1$ and
$B(I_k)$ is a $\delta$-set, $\gamma \in B(I_j)$ implies $\mu \ominus
1 \in B(I_j) \subset B(I_k)$ too, a contradiction. If $d=d_k-1$,
then $\gamma >\mu \ominus 1$. Note that $(\gamma, d_k-1) \in B(I)$
implies $\gamma \in B(I_{k-1})$ by Lemma \ref{d}. However, $\mu \in
\irr(I_{k-1})$ so there is no $\gamma \in B(I_{k-1})$ such that
$\gamma >\mu \ominus 1$, a contradiction. Hence such $(\gamma, d)$
does not exist. Consequently, $(\mu,d_k) \in \irr(I)$.

Conversely, assume $(\mu,d) \in \irr(I)$, we need to prove that there exist some $1\le k \le s$ such that $d=d_k$
and $\mu \in \irr(I_{k-1})\setminus \irr(I_k)$.
By Lemma \ref{irri}, $(\mu,d) \in \irr(I)$ implies
\begin{equation} \label{a}(\mu,d) \ominus 1 \in B(I),
\end{equation}
and there is no $(\gamma,l) \in B(I)$ such that
\begin{equation} \label{b}(\gamma,l)>(\mu,d) \ominus 1.
\end{equation}
By Lemma \ref{d}, (\ref{a}) implies there exists $k$ such that $\mu \ominus 1 \in B(I_{k-1})$, and
\begin{equation}\label{c}
d_{k-1} \le d-1 < d_k.
% \textrm{ ie. }d_{k-1} < d-1 \le d_k.
\end{equation}
By Lemma \ref{d} again, $(\mu \ominus 1,d_k-1) \in \B(I)$. Then
(\ref{b}) and (\ref{c}) imply that $d=d_k$. (\ref{b}) and (\ref{c})
also imply that there is no $\gamma$ such that $\gamma \in
B(I_{k-1})$ and $\gamma >\mu \ominus 1$, so $\mu \in \irr(I_{k-1})$.

It remains to prove $\mu \notin \irr(I_k)$. Assume $\mu \in \irr(I_k)$. Then $\mu \ominus 1 \in B(I_k)$. By Lemma
\ref{d}, $(\mu \ominus 1, d_k) \in \B(I)$ and $(\mu \ominus 1, d_k)>(\mu,d_k) \ominus 1$, contradicting to
$(\mu,d_k) \in \irr(I)$. Thus $\mu \in \irr(I_{k-1})\setminus \irr(I_k)$.
\end{proof}

Theorem \ref{rec} gives us the following recursive algorithm for
finding irreducible decomposition of monomial ideals. Suppose we are
given $I=\langle X^{\alpha_1},\ldots, X^{\alpha_p} \rangle$ and
fixed variable order $x_1 < \cdots <x_n$. We encode the set
$\{\alpha_1, \ldots, \alpha_p\}$ as a tree $\cT$ of height $n$. Our
algorithm $\irr(\cT)$ takes $\cT$ as input and produce $\irr(I)$ as
output. That is, $\irr(I)=\irr(\cT)$.

\vspace{0.2in}
 \noindent{\bf Recursive Algorithm}: $\irr(\cT)$
 \begin{tabbing}
\hspace*{0.55in} \= kk \= xxx \=\kill
Input: \> $\cT$, a tree encoding $\mn(I)$\\
Output: \> $S$, a set (or a tree) representing $\irr(I)$\\
 Step 1. \>  Start at the root of $\cT$. If the height of $\cT$ is $1$, then $\cT$ consists of a few leaves; \\
 \> let $d$ be the largest label on these leaves and let $S:= \{d \}$. \\
 \> Return $S$ (and stop the algorithm).\\
Step 2. \> Now assume $\cT$ has height at least two. Set $S:=\{\ \}$.\\
Step 3. \> Suppose $d_0< d_{1}< \cdots<d_s$ are the labels of the children under the root of $\cT$, \\
\> and let $\cT_k$ be the subtree extending from $d_k$, $0 \le k \le s$. \\
\> Note that the root of $\cT_k$ is the node labeled by $d_k$, but now unlabeled.\\
\> Find $V_{0}:= \irr(\cT_0)$ by recursive call of this algorithm.
\\
\> For $k$ from 1 to $s$ do\\
\> \> 3.1. \> Find $\cT_k:=\minmerge(\cT_{k-1}, \cT_k)$, and delete $\cT_{k-1}$. \\
\> \> 3.2. \> Find $V_{k}:= \irr(\cT_{k})$ by recursive call of this algorithm.\\
\> \> 3.3.  \> Find $V:=V_{k-1}\setminus V_{k}$, delete $V_{k-1}$, and $S:=\merge(S, V \otimes d_{k})$.\\
Step 4. \> Return ($S$).\\
\end{tabbing}

\begin{exmp} \label{3var} We end this section by demonstrating how the algorithm is used to decompose the ideal
$I=\langle x^4,y^4,x^3y^2z^2,xy^3z^2,x^2yz^3 \rangle.$ First
represent the monomials as a tree with variable order $x<y<z$, where
$\cT_k$'s are the subtrees extending from the node with label $d_k$,
$k=0,1,2,3$.
\begin{figure}[h]
$$
\includegraphics*[128,689][565,814]{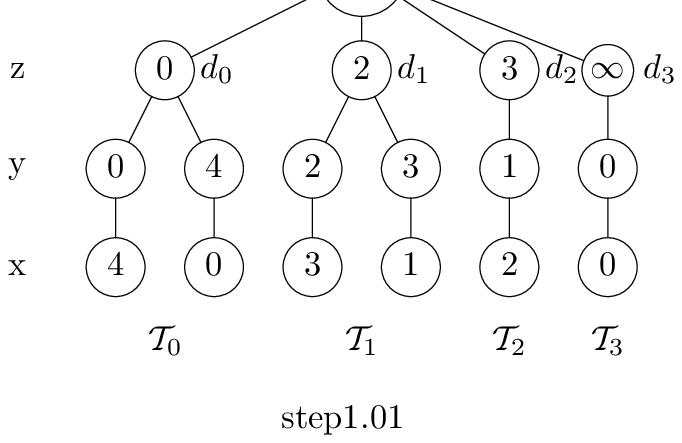}
$$
\caption{Tree representation.}\label{step1}
\end{figure}

\begin{figure}[h]
$$
\includegraphics*[125,689][565,814]{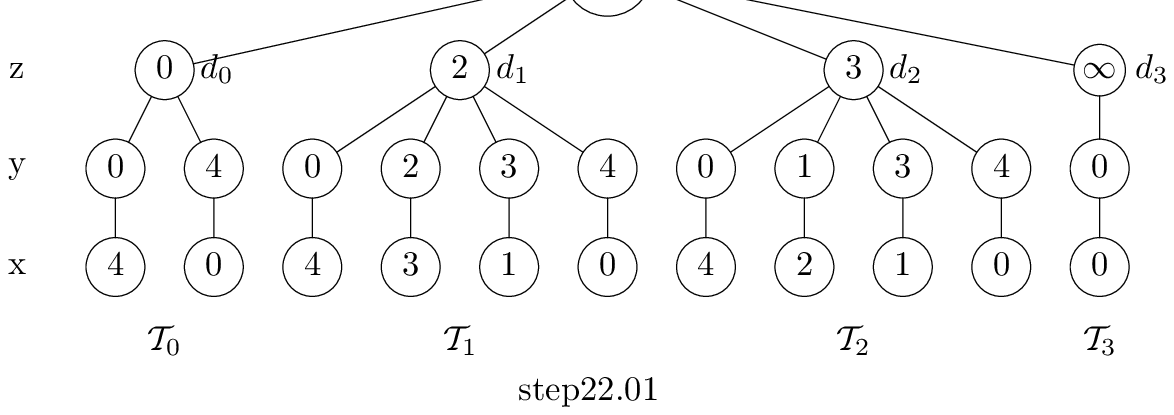}
$$
\caption{MinMerge step.} \label{step2}
\end{figure}
\begin{figure}[h]
$$
\includegraphics*[125,682][565,814]{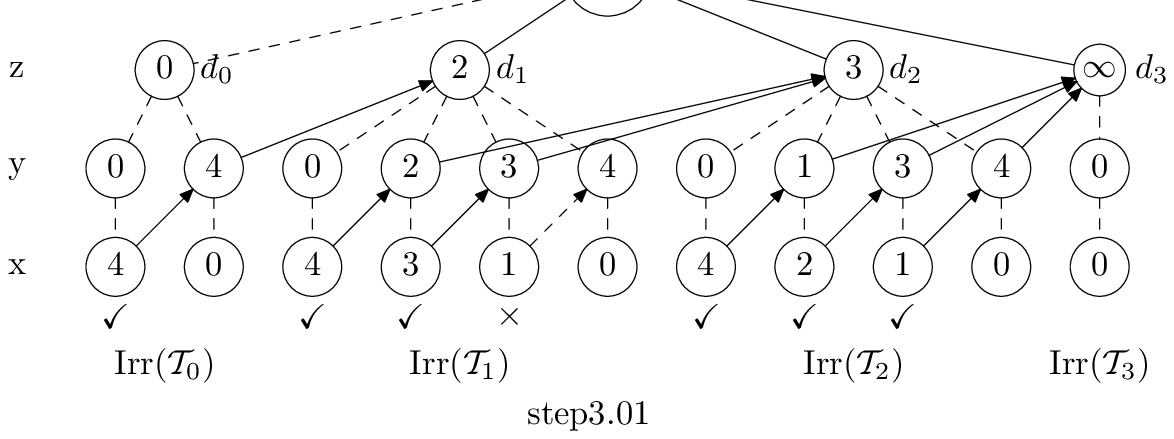}
$$
\caption{Shifting step.}\label{step3}
\end{figure}

Figure \ref{step2}-\ref{step3} show the process of finding the
irredundant irreducible decomposition of $I$. For each $\cT_k$,
inductively MinMerge the subtrees from left to right, corresponding
to Step 3.1 in the Recursive algorithm. See Figure \ref{step2}. In
Figure \ref{step3} we call the procedure $\irr(\ )$ for each $\cT_k$
to compute $\irr(\cT_k)$, corresponding to Step 3.2. Since the
height of $\cT_k$ is 2, we bind each leaf that is not in the
most-right side of $\cT_k$ with the node of height 2 on the next
path - just do the shifting in adjacent paths, see Figure
\ref{step3}. Finally we find the paths in $\irr(\cT_{k-1})$ that are
not in $\irr(\cT_{k})$. The one with a mark $\times$ in
$\irr(\cT_{k})$ is discarded. Then bind the resulting paths with
$d_k$. The irreducible components can be read from the last figure:
$$\irr(I)=\{(4,4,2), (4,2,3), (3,3,3),(4,1,\infty),(2,3,\infty),(1,4,\infty)\}.$$
\end{exmp}

% Distribution ================================================================
\section{Incremental Algorithm} \label{distrbn}
In this section we shall present an incremental algorithm based on
the idea of adding one generator at a time. This algorithm can be
viewed as an improvement of Alexander Dual method
(\cite{Mil00,Mil04}). We maintain an output list of irreducible
components, and at each step we use a new generator to update the
output list. In \cite{Mil04}, it is not clear how to select good
candidates that need to be updated, and the updating process there
is also inefficient. Our algorithm avoids these two deficiencies. We
establish some rules that help us to exclude many unnecessary
comparisons.

Monomial ideal are much simpler than general ideals. The next
theorem tells us that monomial ideals satisfy distribution rules for
the operations \lq\lq$+$\rq\rq\ and \lq\lq$\cap$\rq\rq. These rules
may not be true for general ideals.

\begin{thm}[Distribution Rules] \label{distr} Let $I_1, \ldots, I_t, J$ be any monomial ideals in $\K[X]$.
Then
\begin{enumerate}
\item [(a)] $(I_1+\ldots+I_t)\cap J=I_1\cap J+\ldots+I_t\cap J$, \textrm{and} \item [(b)]$(I_1 \cap \ldots \cap
I_t)+J=(I_1+ J)\cap \ldots \cap (I_t+J)$.
\end{enumerate}
\end{thm}
\begin{proof} By induction, we just need to prove the case for $t=2$. Note that (b) follows form (a), as
\begin{eqnarray*}(I_1+ J)\cap (I_2+J) &=& I_1\cap
(I_2+J)+ J\cap (I_2+J)\\
&=& I_1\cap I_2+I_1\cap J+J\cap I_2+J\\
&=& I_1\cap I_2+J.
\end{eqnarray*}
To prove (a) for the case $t=2$, suppose $h$ is a generator for $(I_1+I_2) \cap J$.
Then $h$ must be in $(I_1+I_2)$ and $J$.
Since $(I_1+I_2)\cap J$ is also a monomial ideal, $h$ is a monomial.
 The fact that $h \in I_{1}+I_{2}$ implies that $h$ is in either $I_{1}$ or $I_{2}$.
 Hence $h$ is in $I_1\cap J$ or in $I_2\cap J$, so $h \in I_1\cap
J+I_2\cap J$. Going backward yields the proof for the other direction. \qedhere
\end{proof}

Theorem \ref{distr} gives us an incremental algorithm for
irreducible decomposition of monomial ideals. Precisely, we have the
following situation at each incremental step: Given the irreducible
decomposition $\irr(I)$ of an arbitrary ideal $I$ and a new monomial
$X^\alpha$ where $\alpha=(a_{1},\ldots,a_{n}) \in \N^n$, we want to
decompose $\widetilde{I}=I+\langle X^{\alpha} \rangle$. By the
distribution rule (b),
\begin{equation} \label{whole}
\widetilde{I} = \left(\bigcap_{\beta\in \irr(I)} m^{\beta}\right) +\langle X^{\alpha} \rangle
 =\bigcap_{\beta \in \irr(I)}\left(m^{\beta}+\langle X^{\alpha}\rangle\right).
\end{equation}
We need to see how to decompose each ideal on the right hand side of (\ref{whole}) and how to get rid of
redundant components. We partition $\irr(I)$ into two disjoint sets:
\begin{eqnarray}
T_1^{\alpha}&=&\{\beta\in \irr(I): \alpha \nprec \beta\}, \label{t11}\textrm{ and}\\
T_2^{\alpha}&=&\{\beta\in \irr(I): \alpha \prec \beta\}.\label{t12}
\end{eqnarray}
Note that if $X^{\alpha} \in I$ then $T_2^{\alpha} =\phi$. For each $\beta \in T_1^{\alpha}$, we have
$X^{\alpha} \in m^{\beta}$, thus
\begin{equation} \label{t1}
m^{\beta}+\langle X^{\alpha}\rangle=m^{\beta}.
\end{equation}
For each $\beta \in T_{2}^{\alpha}$, we have $X^{\alpha} \notin
m^{\beta}$.
 In this case, we split $\langle X^{\alpha} \rangle $ as
$$\langle X^{\alpha}\rangle=\bigcap_{j=1}^n\langle x_j^{a_{j}}\rangle.$$
By the distribution rule (b), we have
$$m^{\beta}+\langle X^{\alpha}\rangle=\bigcap_{j=1}^n \left(m^{\beta}+\langle x_j^{a_{j}}\rangle\right).
$$
Define
$$\beta^{(\alpha,j)}=(b_1,\ldots, b_{j-1},a_{j}, b_{j+1}, \ldots, b_n), \quad 1 \leq j \leq n.$$
Since $\alpha \prec \beta$, we have $a_j < b_j$ for all $1 \le j \le
n$. Hence $m^{\beta}+\langle
x_j^{a_{j}}\rangle=m^{\beta^{(\alpha,j)}}$, and
\begin{equation} \label{t2}
m^{\beta}+\langle X^{\alpha}\rangle=\bigcap_{j=1}^n m^{\beta^{(\alpha,j)}}.
\end{equation}
Therefore,
\begin{equation} \label{wtI}
\irr(\widetilde{I}) = \maxmerge\left(T_1^{\alpha},
\{\beta^{(\alpha,j)}: \beta \in T_2^{\alpha} \textrm{ and } 1\le j
\le n\}\right).
\end{equation}

It remains to see which of the components in the right hand side of
the above expression belong to $\irr(\widetilde{I})$, so others are
redundant.
\begin{lem} \label{case1} $T_1^{\alpha} \subset
\irr(\widetilde{I})$.
\end{lem}
\begin{proof} Let $\beta_1\in T_1^{\alpha}$. By equation (\ref{wtI})
if $\beta_1 \notin \irr(\widetilde{I})$, then there exists some
$\beta_2\in T_2^{\alpha}$ such that $\beta_1$ is maxmergeed by
$\beta_2^{(\alpha,j)}$ for some $j$, ie. $\beta_1 \le
\beta_2^{(\alpha,j)}$. Since $\beta_2^{(\alpha,j)} < \beta_2$,
$\beta_1 \le \beta_2^{(\alpha,j)}$ implies that $\beta_1 < \beta_2$,
which contradicts with the fact that $\beta_1, \beta_2 \in \irr(I)$.
Hence $\beta_1 \in \irr(\widetilde{I})$ as claimed.
\end{proof}
Lemma \ref{case1}  shows that the elements in $T_1^{\alpha}$ will be
automatically in $\irr(\widetilde{I})$. Now we turn to the
components $\beta^{(\alpha,j)}$. For $\beta \in T_2^{\alpha}$,
define
\begin{equation} \label{abeta}
 M_{\beta}=\{m\in \mn(I): m|X^{\beta}\}.
\end{equation}
For $m\in M_{\beta}$, if $\deg_{x_u}m =b_u$, then we say $m$
\textit{matches} $\beta$ in $x_u$. It is possible that one monomial
matches $\beta$ in multiple variables. For example, with $I=\langle
x^2, y^2, z^2, xy, xz, yz \rangle$ and $\beta=(1,1,2) \in \irr(I)$,
the monomial $xy$ matches $\beta$ in $x$ and $y$. We say $m$ matches
 $\beta$ only in $x_u$ if $\deg_{x_u}m =b_u$ and
$\deg_{x_k}m <b_k$ for all $k \neq u$.

\begin{lem} \label{only} For each $\beta=(b_{1}, \ldots, b_{n}) \in T_2^{\alpha}$ and
each $1 \le u \le n$, there exists $m\in M_{\beta}$ such that $m$
matches
 $\beta$ only in $x_u$.
\end{lem}
\begin{proof} Note that a vector $\gamma \in \B(I)$ is maximal if
and only if
 $X^{\gamma} \cdot x_u \in I$ for every $u$.
Since $\beta \in \irr(I)$, $\beta \ominus 1$ is maximal in $B(I)$.
Thus, for each $1 \leq u \leq n$, $X^{\beta \ominus 1} \cdot x_u \in
I$, so there exists a monomial say $m \in \mn(I)$ such that
$m|X^{\beta \ominus 1}\cdot x_u$. Then $\deg_{x_k}m <b_k$ for $k
\neq u$. If $\deg_{x_u}m <b_u$ as well, then $m|X^{\beta \ominus
1}$, which implies that $X^{\beta \ominus 1} \in I$, a
contradiction. Therefore $\deg_{x_u}m =b_u$. Note that $X^{\beta
\ominus 1}\cdot x_u|X^{\beta}$, so $m\in M_{\beta}$.
\end{proof}

For any set of monomials $A \subset \K[X]$, define $\textbf{max}(A)$
be the exponent $\gamma$ such that $X^{\gamma}=\lcm(A)$.

\begin{lem} $\max(M_{\beta})=\beta$.
\end{lem}
\begin{proof} By the definition of $M_{\beta}$,
we know that $\max(M_{\beta})\le \beta$. By Lemma \ref{only} we have
$\max(M_{\beta})\ge \beta$. Thus $\max(M_{\beta})=\beta$.
\end{proof}

For $k\neq u$, let
\begin{equation}\label{dbeta}
d(\beta, u, k) = \min\{\deg_{x_u} m : m \in M_{\beta} \text{
matching $\beta$ only in $x_k$}\}.
\end{equation}
 Note that $d(\beta,u, k) < b_u$. Define
\[ d(\beta,u) = \max_{1\leq k \leq n, k\ne u} \{d(\beta,u,k)\}.\]
\begin{lem} \label{du}
For each $\beta \in T_2^{\alpha}$ and $1 \leq u \leq n$,
$\beta^{(\alpha,u)} \in \irr(\widetilde{I})$ if and only if $d(\beta,u) < a_u$.
\end{lem}
\begin{proof}
Suppose $d(\beta,u) < a_u$. We want to prove that
$\beta^{(\alpha,u)} \in \irr(\widetilde{I})$. By Lemma \ref{irri},
this is equivalent to proving that $\beta^{(\alpha,u)} \ominus 1 \in
\B(\widetilde{I})$ and is maximal. Assume $\beta^{(\alpha,u)}
\ominus 1 \notin \B(\widetilde{I})$. Then there exists $m \in \mn(I)
\cup \{X^{\alpha}\}$ such that $m|X^{\beta^{(\alpha,u)} \ominus 1}$.
First note that $m \neq X^{\alpha}$ because $X^{\alpha}$ can not
divide $X^{\beta^{(\alpha,u)} \ominus 1}$. Thus $m \in \mn(I)$,
which implies $X^{\beta^{(\alpha,u)} \ominus 1} \in I$. Since
$\beta^{(\alpha,u)} \ominus 1<\beta \ominus 1$, we have $X^{\beta
\ominus 1}\in I$, contradicting to $\beta\in \irr(I)$. Hence
$\beta^{(\alpha,u)} \ominus 1 \in \B(\widetilde{I})$. We next need
to prove that $\beta^{(\alpha,u)} \ominus 1$ is maximal in
$\B(\widetilde{I})$, that is, $X^{\beta^{(\alpha,u)} \ominus 1}
\cdot x_k \in \widetilde{I}$ for every $k$. In the case for $k=u$,
we have $X^{\alpha}|X^{\beta^{(\alpha,u)} \ominus 1} \cdot x_u$. For
any $k \neq u$, let $m$ be any monomial in (\ref{dbeta}) such that
$\deg_{x_u} m = d(\beta, u, k)$. Then $\deg_{x_u} m = d(\beta, u, k)
\leq d(\beta,u) < a_u$, hence $m |X^{\beta^{(\alpha,u)} \ominus 1}
\cdot x_k$ as $\deg_{x_k} m =b_k$ and $\deg_{x_j}m \le b_j-1$ for $j
\neq u,k$.

Conversely, suppose $\beta^{(\alpha,u)} \in \irr(\widetilde{I})$. We
want to prove that $d(\beta,u) < a_u$. We know that
$\beta^{(\alpha,u)} \ominus 1$ is maximal in $\B(\widetilde{I})$.
Thus $X^{\beta^{(\alpha,u)} \ominus 1} \cdot x_k \in \widetilde{I}$
for every $k$. For any $k \neq u$, suppose $X^{\beta^{(\alpha,u)}
\ominus 1} \cdot x_k$ is divisible by $m \in \mn(I) \cup
\{X^{\alpha}\}$. Then
\begin{equation} \label{degx_u}
\deg_{x_u}m \le a_u-1 < b_u, \quad \deg_{x_j}m \le b_j-1, \ j \neq
u,k,
\end{equation}
and $\deg_{x_k}m \leq b_k$. As  $X^{\beta^{(\alpha,u)} \ominus 1}
\in \B(\widetilde{I}) \subset \B(I)$, $m$ can not divide
$X^{\beta^{(\alpha,u)} \ominus 1}$. Hence $\deg_{x_k}m \leq b_k$. So
$m$ matches $\beta$ only in $x_k$. Note that $m \neq X^{\alpha}$, so
$m \in M$ and thus $m \in M_{\beta}$. It follows that $d(\beta,u,k)
\le a_u-1$ by (\ref{degx_u}). Therefore, $d(\beta,u) < a_u$ as
desired.
\end{proof}

By the above lemma, for each $\beta \in T_2^{\alpha}$, we only need
to find $M_{\beta}$ and $d(\beta,u)$, which will tell us whether
$\beta^{(\alpha,u)} \in \irr(\widetilde{I})$. This gives us the following
incremental algorithm.

\vspace{0.2in} \noindent {\bf Incremental algorithm}
\begin{tabbing}
\hspace*{0.55in} \= kk \= kkk \= \hspace*{0.5in} \=\kill
Input: \> $M$, a set of monomials in $n$ variables $x_1,\ldots, x_n$.\\
Output: \> $\irr(I)$, the irredundant irreducible components of the
ideal $I$ generated by $M$.\\
 Step 1. \>  Compute $\minmerge(M)$ and sort it into the form: \\
 \> \> \> \> $\minmerge(M)=\{x_1^{c_1},\ldots, x_n^{c_n},X^{\alpha_1},\ldots,
X^{\alpha_p}\},$\\
 \> where $c_i$ can be $\infty$ and
$\{X^{\alpha_1},\ldots, X^{\alpha_p}\}$ are sorted in lex order with variable\\
 \> order $x_1<\ldots<x_n$. Set\\
 \>  \> \> \> $T:=\{(a_1, \ldots, a_n)\}$.\\
Step 2.  \> For each $k$ from 1 to $p$ do:\\
 \> \> 2.1. Set the temporal variables $V=\emptyset$ and $\alpha :=\alpha_k$.\\
 \> \> 2.2. For every $\beta \in T$ with $\alpha \nprec \beta$ do\\
 \>  \> \> \> $V:=V \cup \{\beta\}.$\\
 \> \> 2.3. For every $\beta \in T$ with $\alpha \prec \beta$ do,\\
 \> \> \> $\bullet$ find $M_{\beta}$ as defined in (\ref{abeta});\\
 \> \> \> $\bullet$ for $1 \le u \le n$, compute $d(\beta, u)$, and if $d(\beta,u) < a_u$ then update\\
 \>  \> \> \> $V:=V \cup \{\beta^{(\alpha,u)}\}.$\\
 \> \> 2.4. Set $T:=V$.\\
Step 3. \> Output $T$.\\
\end{tabbing}

We next prove that there is a nice property of the above algorithm for
generic monomial ideals, that is, the size of $T$ is always non-decreasing
at each stage when a new generator is added. This will allow us to bound
the running time of the algorithm in term of input and output sizes.

\begin{thm} \label{nondecr} Suppose $I$ is generic and
 $\mn(I)=\{ x_1^{c_1}, \ldots, x_n^{c_n},
X^{\alpha_1},\ldots,X^{\alpha_p}\}$ where $X^{\alpha_k}$'s are
sorted in lex order with variable order $x_1<\ldots<x_n$. Let
$\widehat{I}=\langle x_1^{c_1}, \ldots, x_n^{c_n},$
$X^{\alpha_1},$$\ldots,X^{\alpha_{p-1}} \rangle$. Then
$|\irr(\widehat{I})| \le |\irr(I)|$.
\end{thm}
\begin{proof} Keep notations as above.
For every $\beta \in T_2^\alpha$, $b_n=c_n$. Thus $x_n^{c_n}$ is the
only monomial in $M_\beta$ that has degree in $x_n$ larger than
$a_n$. Hence $d(\beta, n) <a_n$ and $\beta^{(\alpha,n)} \in
\irr(I)$. By the equation (\ref{wtI}) and Lemma \ref{case1},
$$|\irr(I)| \ge |T_1^{\alpha}|+|\{\beta^{(\alpha,n)}: \beta \in T_2^{\alpha}\}|=|T_1^{\alpha}|+|T_2^{\alpha}|= |\irr(\widehat{I})|.\qedhere$$
\end{proof}

The reader might wonder
whether a similar statement holds in non-generic case as well.
 The answer is negative. Let $I=\langle
x^3,y^3,z^2,w^2,x^2yz,xy^2w\rangle \subset \K[x,y,z,w]$ with lex
order and $x<y<z<w$. Then
$$\irr(I)=\{(3,3,1,1),(2,3,2,1),(3,2,1,2),(3,1,2,2),(2,2,2,2),(1,3,2,2)\}.$$
By adding $X^\alpha=xyzw$, we can see $\beta =(2,2,2,2)\in
T_2^\alpha$. Note that $M_\beta=\{x^2yz,xy^2w,z^2,$ $w^2\}$. Since
$d(\beta, u) =1=a_u$ for $u=1,2,3,4$, no new $\beta^{(\alpha,j)}$
will be generated. Thus the number of irreducible components
decreases by 1 instead.

We find the irreducible components for the monomial ideal in Example
\ref{3var} again by the flow of our incremental algorithm.

\begin{exmp} \label{3var2} Decompose $$I=\langle x^4,y^4,x^3y^2z^2,xy^3z^2,x^2yz^3 \rangle.$$

Note: \lq\lq$\checkmark$\rq\rq\ means $\beta^{(\alpha,u)} \in
\irr(\widetilde{I})$ for corresponding $\beta, \alpha$ and $u$,
while \lq\lq$\times $\rq\rq\ means not.

\begin{tabbing}
\hspace*{0.5in} \= kk \= kkk \= \hspace*{0.08in} \=
\hspace*{0.04in}\=\kill
Step 1.  \> $M=\{x^4,y^4,z^{\infty}, x^3y^2z^2,xy^3z^2,x^2yz^3\}$. Set $T:=\{(4, 4, \infty)\}.$\\
Step 2.  \> (i) For $\alpha=(3,2,2)$ do:\\
 \> \> 2.1. $V:=\phi$.\\
 \> \> 2.2.  Since $\alpha \prec (4, 4, \infty)$, $V:=\phi$.\\
 \> \> 2.3.  Let $\beta=(4, 4, \infty)$. We find $M_{\beta}=\{x^4,y^4\}$.\\
 \>  \> \>So we have $d\{\beta,1\}=0$($\checkmark$), $d\{\beta,2\}=0$($\checkmark$) and
$d\{\beta,3\}=0$($\checkmark$).\\
 \>  \> \>Then $V:=\{(3,4,\infty),(4,2,\infty),(4,4,2)\}.$\\
 \> \> 2.4.  Let $T:=V$.\\
\> (ii) For $\alpha=(1,3,2)$ do:\\
 \> \> 2.1.  $V:=\phi$.\\
 \> \> 2.2.  Update $V$ by $V:=\{(4,4,2),(4,2,\infty)\}$.\\
 \> \> 2.3.  $\alpha \prec (3,4,\infty)$.\\
 \>  \> \>Let $\beta=(3,4,\infty)$. We find $M_{\beta}=\{y^4,x^3y^2z^2\}$.\\
 \>  \> \>So $d\{\beta,1\}=0$($\checkmark$), $d\{\beta,2\}=2$($\checkmark$)
and $d\{\beta,3\}=2$($\times$). \\
 \>  \> \>Then
 $V:=\{(4,4,2),(4,2,\infty),(1,4,\infty),(3,3,\infty)\}.$\\
 \> \> 2.4.  Let $T:=V$.\\
\> (iii) For $\alpha=(2,1,3)$ do:\\
 \> \> 2.1.  $V:=\phi$.\\
 \> \> 2.2.  $V:=\{(4,4,2),(1,4,\infty)\}$.\\
 \> \> 2.3.  $\alpha \prec
(4,2,\infty)$, and $\alpha \prec(3,3,\infty)$.\\
 \> \> \> $\bullet$ Let $\beta=(4,2,\infty)$. We find
$M_{\beta}=\{x^4,x^3y^2z^2\}$.\\
 \>  \> \> \>So $d\{\beta,1\}=3$($\times$), $d\{\beta,2\}=0$($\checkmark$) and
$d\{\beta,3\}=2$($\checkmark$). \\
 \>  \> \> \>Then
$V:=\{(4,4,2),(1,4,\infty),(4,1,\infty),(4,2,3)\}.$\\
 \> \> \> $\bullet$  Let $\beta=(3,3,\infty)$. Then $M_{\beta}=\{x^3y^2z^2,xy^3z^2\}$.\\
 \>  \> \> \>$d\{\beta,1\}=1$($\checkmark$), $d\{\beta,2\}=2$($\times$),
$d\{\beta,3\}=2$($\checkmark$). \\
 \>  \> \> \>So
$V:=\{(4,4,2),(1,4,\infty),(4,1,\infty),(4,2,3),(2,3,\infty),(3,3,3)\}.$\\
 \> \> 2.4.  Let $T:=V$.\\
Step 3. \> Output $T$ \> \> \>\>$=\{(4,4,2),(1,4,\infty),(4,1,\infty),(4,2,3),(2,3,\infty),(3,3,3)\}$\\
 \>  \> \> \>\>$=\{(4,4,2), (4,2,3),
(3,3,3),(4,1,\infty),(2,3,\infty),(1,4,\infty)\}$.\\
\end{tabbing}
\end{exmp}

Some preprocess can be taken right before Step 2 to improve the
efficiency of the incremental algorithm. For each $u \in \{1,
\ldots, n\}$, we partition $M$ into disjoint subsets such that the
monomials in each subset have the same degree in $x_u$. We then
store these information, which requires memory complexity $O(n \cdot
p)$. For each $\beta \in T_2^{\alpha}$, we can find $M_{\beta}$ by
only checking the monomials in the subset with degree $b_u$ in
variable $x_u$ for every $u$. Note that for generic monomial ideals
each subset contains a unique monomial. In this case $M_{\beta}$
contains $n$ monomials, and it can be found by $O(n)$ operations,
instead of
$O(p)$ operations by scanning through the whole input monomial set.\\

\section{Time Complexity and Conclusion} \label{Conclusion}

We estimate the running time of our algorithms by counting the
number of monomial operations (ie. comparisons and divisibility)
used. Our recursive algorithm depends heavily on the number of
distinct degrees of each variable. Let $s_j$ be the number of
distinct degrees of $x_j$ where $j=1,\ldots,n$. Then the total
number of merge of subtrees used by the algorithm is at most
$\prod_{j=1}^n s_j$. Since each subtree has at most $p$ leaves(ie.
$p$ generators), each merge takes $O(p^2)$ monomial operations.
Hence the algorithm uses $O(p^2 \cdot \prod_{j=1}^n s_j)$ monomial
operations. This algorithm is more efficient for highly non-generic
monomial ideals. The benchmark analysis in \cite{Rou08b} compare
several algorithms based on various slicing strategies, including
our recursive algorithm. It is shown there that our algorithm
performs as a very close second best one.

The running time of our incremental algorithm is harder to estimate
for general ideals. For generic ideals, however, we can bound the time
in terms of input and output sizes. More precisely,
suppose
\[ I=\langle x_1^{c_1}, \ldots, x_n^{c_n}, X^{\alpha_1},\ldots,X^{\alpha_p}\rangle \]
is a generic monomial ideal in $\K[X]$ where $X^{\alpha_k}$'s are
sorted in lex order with variable order $x_1<\ldots<x_n$. For $0
\leq k \leq p$, let
\[I_{(k)} = \langle x_1^{c_1}, \ldots, x_n^{c_n}, X^{\alpha_1},\ldots,X^{\alpha_k}\rangle. \]
All these ideals are generic. By Theorem \ref{nondecr},
we have
\[ 1=|\irr(I_{(0)})| \leq |\irr(I_{(1)})| \leq \cdots \leq |\irr(I_{(p)})| = |\irr(I)|.\]
In an arbitrary stage of the incremental algorithm, we try to find
the irreducible components of $I_{(k)}$ from those of $I_{(k-1)}$.
For each $\beta \in \irr(I_{(k-1)})$, only those $\beta$ in
$T_2^{\alpha_k}$ (as defined in (\ref{t12})) need to be updated.
Note that $I$ is generic, by the preprocess $M_{\beta}$ can be found
in $O(n)$ operations. The numbers $d(\beta, u, k)$, $1\leq u \ne k
\leq n$, can be computed by scanning through the monomials in
$M_{\beta}$ once, thus using only $O(n)$ monomial operations. Then
the numbers $d(\beta, u)$, $1\leq u \leq n$, can be computed in
$O(n^2)$ operations. Hence for each $\beta \in T_2^{\alpha_k}$, Step
2.3 uses at most $O(n+n^2)=O(n^2)$ monomial operations. Since $T
\supset T_2^{\alpha_k}$ has at most $\ell$ elements where
$\ell=|\irr(I)|$, Step 2.3 needs at most $O(n^2\ell)$ monomial
operations. Therefore, the total number of monomial operations is at
most $O(n^2p\ell)$. In fact, $T_2^{\alpha_k}$ is usually a small
subset of $T$, the actual running time is much better than our
worst-case estimate indicates.

We also want to point out that for generic monomial ideals, the
incremental algorithm is an improved version of the recursive
algorithm. Suppose we add the new monomial $X^{\alpha_k}$ into
$I_{(k-1)}$. In Step 3.2 of the recursive algorithm, we need to
compute $\irr(\cT_k)$. But in Step 2.3 of the incremental algorithm,
only $\beta \in T_2^{\alpha_k}$ need to be updated. We have the
observation that $T_2^{\alpha_k}$ is a small subset of
$\irr(\cT_k)\otimes c_n$. By this observation we conclude the
incremental algorithm is more efficient than the recursive algorithm
for generic monomial ideals. In non-generic case, the comparison is
not clear.

In all previous algorithms (including our recursive one) for
monomial decomposition, the storage in the intermediate stages may
grow exponentially larger than the output size. Our incremental
algorithm seems to be the first algorithm for monomial decomposition
that the intermediate storage is bounded by the final output size.
Note that the output size $\ell$ can be exponentially large in $n$.
In fact, it is proven in \cite{Agn97} that $\ell =
O(p^{[\frac{n}{2}]})$ for large $p$. Since the output size can be
exponential in $n$, it is impossible to have a polynomial time
algorithm for monomial decomposition.\\

\section{Acknowledgement}

We thank Alexander Milowski and Bjarke Roune for comments and
suggestions, and Ezara Miller for helpful communications (especially
for providing some of the diagrams).

\end{document}